\documentclass[12pt]{article}
\usepackage{ctable}
\usepackage{multirow}
\usepackage{subfigure}
\usepackage{amssymb}
\usepackage{amsmath}
\usepackage{amsthm}
\usepackage{rotating}
\usepackage{natbib}

\DeclareMathOperator*{\E}{\mathbb{E}}
\DeclareMathOperator*{\sto}{\operatorname{subject}\,\operatorname{to}}

\title{A comparison of non-stationary stochastic lot-sizing strategies}
	
	\author{Gozdem Dural-Selcuk\\ 
	{\scriptsize Institute of Population Studies, Hacettepe University, Ankara, Turkey}\\
	\and
	Onur A. Kilic\\
	{\scriptsize Institute of Population Studies, Hacettepe University, Ankara, Turkey}\\
	\and
	S. Armagan Tarim\\
	{\scriptsize Institute of Population Studies, Hacettepe University, Ankara, Turkey}\\
	\and
	Roberto Rossi\\
	{\scriptsize Business School, University of Edinburgh, Edinburgh, United Kingdom}
	}

\begin{document}
	
	\maketitle
	
	\begin{abstract}
	We consider the non-stationary stochastic lot sizing problem with backorder costs and make a cost comparison among different lot-sizing strategies. We initially provide an overview of the strategies and some corresponding solution approaches in the literature. We then compare the cost performances of the lot-sizing strategies on a common test bed while taking into account the added value of realized demand information. The results of this numerical experience enable us to derive novel insights about the cost performance of different stochastic lot-sizing strategies under re-planning with respect to demand realization.
	\end{abstract}

\section{Introduction}
\label{sec:introduction}

The lot sizing problem aims to determine when and how much to order so as to match supply and demand while minimizing inventory costs \citep{Silver1998}. A stream within the literature on the lot sizing problem addresses non-stationary stochastic demand. The assumption of non-stationarity stems from short product life cycles: as product life cycles get shorter, demand rates are subject to quick changes over time. Furthermore, due to seasonality and trend, demand rates may even change within the phases of the product life cycle. As a result, demand is usually not only stochastic but also non-stationary; that is, its probability distribution varies over time \citep{Graves2008}. When demand is stochastic and non-stationary the decision maker must employ inventory control policies with time-varying parameters in order to hedge against demand uncertainty. As shown in \citep{Tunc2011}, it is costly to adopt stationary policies when demand is non-stationary.

Our study focuses on the non-stationary stochastic lot sizing problem under complete backordering and penalty cost. This is a well-studied problem in the literature. In his seminal work \citep{Scarf1960}, Scarf demonstrated the optimality of ($s$,$S$)-type policies under this problem setting. The computation of optimal stationary ($s$,$S$) policy parameters was investigated shortly thereafter in \cite{Veinott1965}. However, although the form of the optimal policy has been known for a long time, computing optimal non-stationary ($s$,$S$) policy parameters still remains a computationally challenging task. Early heuristics that were proposed to address this problem include \cite{Askin1981}, which represents an extension of the classical Silver-Meal \citep{Silver1973} algorithm to a stochastic setting, \cite{Bookbinder1988} and \cite{Bollapragada1999}. The framework introduced by \cite{Bookbinder1988}, which comprises three types of policies, static-uncertainty, static-dynamic uncertainty, and dynamic uncertainty, motivated a number of follow up works, such as \cite{Sox1997,Vargas2009}, who focused on the static-uncertainty policy; and \cite{Tarim2006,Rossi2012}, who focused on the static-dynamic uncertainty strategy. In addition to those mentioned, many other works appeared in the literature, which tackled variants of the problem such as alpha service level constraints \citep{Tarim2004,Tempelmeier2007,Rossi2008,Tarim2011}, beta (i.e. fill rate) service level constraints \citep{Tempelmeier2011}, probabilististic delivery lead time \citep{citeulike:7288002}, etc. All these variants are however out of the scope of our work.

In an early work, \cite{Bookbinder1986} contrasted the performance of six separate lot-sizing methods for use in a rolling schedule. To the best of our knowledge, there is no a similar study focusing on state-of-the-art stochastic lot sizing heuristics under complete backordering and penalty cost. Moreover, the deployment of stochastic lot sizing heuristics in a rolling horizon and/or re-planning setting has yet to be investigated.
 
The current manuscript aims at filling some of these gaps by making the following contributions to the literature:
\begin{itemize}
\item We compare state-of-the-art solution approaches for the non-stationary stochastic lot sizing problem with backorder cost/penalty cost over a common test bed;
\item For the first time in the literature we contrast for the performance of these solution approaches under re-planning; 
\item We derive a number of managerial insights; in particular we demonstrate that, under re-planning, Bookbinder and Tan's static-uncertainty policy, which is known to perform poorly compared to other existing policies, turns out to be near-optimal.
\end{itemize}

The rest of the paper is structured as follows. Section~\ref{sec:background} gives a formal definition of the problem under consideration. Section \ref{sec:survey} gives a detailed overview of the existing approaches in the literature. Section~\ref{sec:numerical} presents the numerical study and the results obtained. Finally, Section~\ref{sec:conclusions} concludes with final remarks.

\section{Background and Literature Review}
\label{sec:background}

We focus on the single-item single-stocking location stochastic lot sizing problem, which is defined as follows. The planning horizon is finite and it consists of $N$ discrete time periods. Individual period demands $d_1,d_2,\ldots,d_N$ are independent random variables with known probability distribution functions not necessarily identically distributed.  A holding cost, $h$ is incurred on any unit carried in inventory from one period to the next. Any demand that cannot be satisfied immediately is backordered and satisfied whenever a sufficiently sized replenishment order arrives. A penalty cost $b$ is incurred for each unit of demand backordered per period. A fixed ordering cost $K$ is incurred every time an order is placed. 

A mathematical formulation for the problem may be written as follows \citep[see e.g.][]{Sox1997,Tarim2006}:
\begin{align*}
	\min \quad & \sum_{n=1}^{N} \E\left\{K z_n + h x_n^+ + b x_n^-\right\} \\
	\sto \quad 
	& x_n = x_{n-1} + q_n - d_n \quad \forall n \in [1,N] \\
	& z_n = \begin{cases} 1 & \text{if $q_n > 0$} \\ 0 & \text{otherwise} \end{cases} \quad \forall n \in [1,N] \\
	& z_n \in \{0,1\},\; q_n \in \mathbb{R}_+ \quad \forall n \in [1,N]
\end{align*}
Here $z_n$ and $q_n$ are the decision variables indicating the replenishment action and the replenishment quantity for period $n$; where $x_n$ stands for the end of period inventory position. Also, $x^+=\max(0,x)$ and $x^-=\max(0,-x)$ respectively stand for the amount of excess stocks and backorders.

\citet{Bookbinder1988} argue that the aforementioned model leads to different control strategies depending on when replenishment decisions $z_n$ and $q_n$ are made. They discuss three classes of strategies: dynamic uncertainty, static uncertainty, and static-dynamic uncertainty. In the dynamic uncertainty strategy, the decision maker observes the current inventory position at the beginning of each time period $n$, i.e. $x_{n-1}$, and decides whether to place an order, and if so how much to order. Therefore, $z_n$ and $q_n$ are determined at period $n$. In the static uncertainty strategy, timing and quantity of orders are fixed once and for all at the beginning of the planning horizon. Thus, all $q_n$ and $z_n$ are determined at the outset. In the static-dynamic uncertainty strategy, the complete replenishment schedule is determined at the beginning of the planning horizon, whereas order quantities are decided at the time of replenishments. Hence, all $z_n$ are determined at the outset but the decision on $q_n$ is postponed until period $n$. 

Intuitively, the cost-effectiveness of an inventory control policy is positively correlated with the amount of demand information used at the time of making replenishment decisions. Thus, the dynamic uncertainty strategy has a superior cost performance as compared to static-dynamic and static uncertainty strategies since it makes replenishment decisions only after observing the actual inventory level. \citet{Tunc2013} numerically analyse the cost performance of these strategies. They show that the static-dynamic uncertainty strategy is very competitive, whereas the static uncertainty performs rather poor as compared to the dynamic uncertainty strategy. As we will show in our numerical study, results are substantially different under a re-planning deployment of these policies.

\section{Survey of Exact and Heuristics Approaches in the Literature}
\label{sec:survey}

Next, we review both optimal and heuristic approaches proposed in the literature to compute the parameters of dynamic, static and static-dynamic uncertainty strategies for the stochastic lot sizing problem with and without re-planning approach. For the purposes of this study, we specifically address those studies analysing non-stationary stochastic inventory systems operating under a penalty cost scheme for backordered demand. So, the solution approaches subject to our analysis is outlined by this exact problem definition.   

\subsection{Dynamic Uncertainty Strategy} 

The literature on dynamic demand dates back to \citeauthor{Wagner1958}'s (\citeyear{Wagner1958}) well-known deterministic dynamic lot sizing model. In his seminal paper, \citet{Scarf1960} characterizes the optimal control policy for the lot sizing problem under stochastic demand; this characterization holds for stationary as well as non-stationary demand. The optimal policy determines two critical parameters for each period: the re-order level and the order up-to level. The decision maker observes the inventory position at the beginning of a period, and places an order if it is below the re-order level so as to replenish the inventory up to the order up-to level. If the inventory position is above the reorder level no order is placed. The optimal policy follows the dynamic uncertainty strategy because both the timing and the quantity of orders become known only at replenishment epochs. 

\citet{Scarf1960} identified the structure of the optimal policy, but finding optimal parameters of the optimal policy remained a computationally intensive task. That is because one needs to recursively compute a continuous cost function in order to obtain the optimal re-order and order-up-to levels for each and every period within the planning horizon. An alternative approach that can be used to tackle the continuity issue is to use a discrete demand distribution \citep[see e.g.][]{Bollapragada1999}. In this case, it is possible to use a discrete state space dynamic program to obtain the optimal cost function. Nevertheless, the resulting procedure is still complex since there is a very wide range of possible inventory levels that should be considered. A variety of studies address this issue \citep[see e.g.][]{Veinott1965,Ehrhardt1979,Federgruen1984,Zheng1991}. However, most of these handle the problem under stationary demand and an infinite planning horizon. 

There are a few studies suggesting heuristic methods to compute parameters of the dynamic uncertainty strategy. \citet{Askin1981} proposes a heuristic method which adapts \citeauthor{Silver1973}'s (\citeyear{Silver1973}) well-known heuristic designed for the deterministic version of the same problem. \citeauthor{Askin1981}'s heuristic first determines the cycle length, i.e. the number of periods to be covered, and the order-up-to level for all prospective orders through the planning horizon. The length of a replenishment cycle is selected so as to minimize the corresponding average total cost per period. The re-order level associated with a replenishment cycle, on the other hand, is myopically set to an inventory level which minimizes the costs to be incurred during the imminent replenishment cycle. The heuristic determines the re-order levels by means of a trade-off analysis between expected costs per period in cases of ordering and not ordering. In particular, the re-order level of a particular period is set to an inventory level where the difference between the expected costs of ordering and not ordering equals the fixed ordering cost. \citet{Bollapragada1999}, on the other hand, propose a myopic heuristic to compute the parameters of dynamic uncertainty strategy. This heuristic relies on the idea of approximating the non-stationary problem by a series of stationary problems. The heuristic proceeds as follows. First, by means of the method developed by \citet{Zheng1991}, for all possible values of mean demand, the optimal parameters of the associated stationary problem as well as the expected time between two consecutive orders are obtained and tabulated. Then, for each period, the order-up-to and re-order levels are set to the corresponding optimal parameters of a specific stationary problem which is chosen in such a way that the cumulative mean demands of the stationary and non-stationary problems over the expected reorder cycle of the stationary problem are equal to each other. The relevant values used in this procedure are read from the table generated in the first step of the heuristic. An important drawback of this method is that it cannot account for end-of-horizon effects due to the underlying stationary approximation. \citeauthor{Bollapragada1999} overcome this issue by replacing the heuristic parameters with the optimal ones for some periods at the end of the planning horizon.

\subsection{Static Uncertainty Strategy}

The static uncertainty strategy makes timing and quantity decisions regarding replenishments at the beginning of the planning horizon. It is customary in the literature to denote the timing as $R$ and the quantity as $Q$; for this reason this policy is often denoted as the ($R$,$Q$) policy. These decisions are not revised in response to realized demands over the planning horizon, and the demand uncertainty is rather absorbed through dimensioning of the fixed lot sizes. However, the static uncertainty strategy has the advantage of providing a completely stable production environment. It is, therefore, appealing in industrial environments characterized by a low degree of flexibility. 

 \citet{Sox1997} studies this strategy and initially models the problem as a mixed integer non-linear program. Then, he provides a solution algorithm based on a network formulation of the problem where each arc corresponds to a prospective replenishment cycle. Here, arc costs are calculated by using optimal cumulative replenishment quantities up to and including the corresponding replenishment cycles. \citeauthor{Sox1997} also derives some properties of the optimal solution which he uses to increase the efficiency of the proposed algorithm. \citeauthor{Sox1997}'s algorithm can be regarded as a stochastic extension of the \citeauthor{Wagner1958}'s (\citeyear{Wagner1958})  algorithm with additional feasibility constraints. \citet{Vargas2009} employs the same approach and shows that optimal replenishment quantities for a given replenishment schedule follow a critical ratio rule. Also, he explicitly addresses the special case where demands are normally distributed, and provide a simple, yet very efficient solution algorithm by exploiting the properties of the normal distribution. 

\subsection{Static-Dynamic Uncertainty Strategy}

The static-dynamic uncertainty strategy provides a stable replenishment pattern by fixing the timing of future orders in advance \citep{Tarim2006}; however, actual order quantities are decided only after demand during periods that precede a given replenishment have been realised. To do so, the decision maker fixes a so-called ``order-up-to-level'' for each replenishment, this represents the level up to which inventory must be raised by the current order. It is customary in the literature to denote the order-up-to-level as $S$; for this reason this policy is often denoted as the ($R$,$S$) policy. This policy is less conservative as compared to the static uncertainty strategy, and it can hedge against uncertainty more efficiently \citep{Kilic2011,Tunc2013}. Because it offers a fixed order schedule, the static-dynamic strategy is particularly appealing in material requirement planning, joint replenishment, and shipment consolidation environments \citep[see][]{Kingsman1985,Silver1998,Li2009,Mutlu2010}. 

The contributions in this line of research are as follows. \citet{Silver1978} proposes a heuristic which can be regarded as the stochastic version of the \citeauthor{Silver1973} (\citeyear{Silver1973}). As is the case in the \citeauthor{Silver1973}'s approach, this heuristic sequentially determines the timing of the next replenishment period starting from the first period of the planning horizon. Here, the expected cost per period is defined as a function of the number of periods the current order is to cover where the associated order-up-to level is myopically determined so as to minimize the expected costs to be incurred until the next replenishment epoch. The procedure postpones the next replenishment period as long as the expected cost per period is decreasing. In this procedure, penalty costs for backordered demands are not explicitly mentioned and order quantities are rather determined to ensure a desired service level until the next replenishment epoch. The service level implementation makes this solution approach to lie beyond the scope of our numerical analysis that especially remarks the non-stationary stochastic lot-sizing problem under penalty cost scheme. Likewise, in what follows we will not survey works that deviate from the assumptions of complete demand backordering and penalty cost.

\citet{Tarim2006} provide a mixed integer programming (MIP) formulation of the problem. As opposed to the procedure employed by \citet{Silver1978}, their method determines order-periods and order-up-to levels simultaneously under the penalty cost assumption. This, however, comes with an additional difficulty which stems from the interdependence between order-up-to levels and costs associated with consecutive replenishment cycles, as well as the non-linear cost function. \citeauthor{Tarim2006} overcome these issues by making the additional assumption that demands are normally distributed. This enables them to develop a certainty equivalent mixed integer programming model where a piecewise linear approximation is used to express the non-linear terms in the objective function. \citeauthor{Tarim2006}'s formulation does not allow expected order sizes to be negative. Hence, if the inventory level at the beginning of a replenishment cycle happens to exceed the corresponding order-up-to level, then the excess inventory is carried forward. Nevertheless, their method implicitly assumes that such instances are rare events and ignores the cost of carrying excess inventories. As such, their approach does not necessarily provide the optimal solution. \citet{Rossi2012} develop a stochastic constraint programming model which is mainly equivalent to \citeauthor{Tarim2006}'s (\citeyear{Tarim2006}) MIP model. However, instead of employing a piecewise linear approximation, they use the exact non-linear cost function. They also introduce a cost based filtering method \citep[see e.g.][]{Focacci2002,Rossi2008a} which exploits the convexity of the cost-function for a fixed replenishment schedule. The filtering method dynamically produces bounds on the optimal total cost during the search procedure of the constraint program for fixed values of binary variables indicating the timing of replenishment epochs, and hence, leads to significant improvements in the computational performance of the proposed method. Although \citeauthor{Rossi2012}'s formulation makes use of the exact non-linear cost function, it may not always yield the optimal policy parameters because -- as is the case for \citeauthor{Tarim2006}'s formulation -- it ignores the cost of carrying excess inventories in cases where the actual inventory level exceeds the order-up-to level. More recently, \cite{citeulike:13341691} extended \citeauthor{Tarim2006}'s MIP approach to generic demand distributions and to a number of service level measures.

\subsection{Computational efficiency of existing heuristics}

We shall now briefly focus on the computational complexity of the approaches considered in our study. Table \ref{tab:overview} lists all approaches considered in our study as well as their computational complexity.
\begin{table}[htbp]
	\centering
	\small
	\begin{tabular}{l|l|l|l|l}
		\textbf{Author(s)} & \textbf{Abbrv.} & \textbf{Strategy} & \textbf{Approach} & \textbf{Complexity} \\
		\toprule
		\citet{Askin1981} & Ask & Dynamic & Heuristic & O($N$) \\
		\citet{Bollapragada1999} & Bol & Dynamic & Heuristic & O($N$) + O($M$) \\
		\citet{Tarim2006} & Tar & Static-Dynamic & Heuristic & Comb. (MIP)\\
		\citet{Rossi2012} & Ros & Static-Dynamic & Heuristic & Comb. (MIP) \\
		\cite{Sox1997} & Sox/Var & Static & Optimal & O($N^2$) (DP) \\
		\cite{Vargas2009} & Sox/Var & Static & Optimal & Comb. (MIP)\\
		\bottomrule
	\end{tabular}
	\caption{An overview of the methods considered in the current study; ``Comb.'' denotes a combinatorial complexity.}
	\label{tab:overview}
\end{table}
As we can observe, with a linear complexity \citeauthor{Askin1981}'s approach is computationally attractive; the approach processes all $N$ periods in the planning horizon sequentially and utilises a closed form critical fractile expression --- similar to a newsvendor solution --- to determine the optimal order quantity. \citeauthor{Bollapragada1999}'s approach is also very efficient: the heuristic essentially boils down to a line search procedure over a pre-computed table. However, the precomputation of the table has pseudo-polynomial complexity, since table elements should cover all possible values that the expected demand may take at a given period up to the maximum expected demand $M$. \citet{Tarim2006}, \citet{Rossi2012} and \citet{Vargas2009} rely on MIP formulations; while \cite{Sox1997} proposed a heuristic with quadratic time complexity that takes the form of a forward dynamic programming (DP) algorithm. 

All these approaches, including those based on MIP formulations, are computationally very efficient and can solve realistic instances with up to 30 periods in fractions of a second, as illustrated in the studies listed in Table \ref{tab:overview} as well as in related follow up works. Conversely, an exact stochastic dynamic programming approach has pseudo-polynomial complexity and, depending on the demand and state space discretization/truncation adopted, it may result computationally very expensive. Since differences in computational performance of all approaches listed in Table \ref{tab:overview} are negligible for most practical purposes, in what follows we will only focus on their relative cost performance.

\section{Numerical Study}
\label{sec:numerical}

This numerical study aims to present a cost based comparison of stochastic lot-sizing strategies under non-stationary stochastic demand and penalty cost assumptions on a common test bed from the literature. The list of the approaches in the scope of this numerical experiment is given in Table~\ref{tab:overview}.  

In addition to comparing the approaches as such, we also compare them under a re-planning deployment, which 
is widely used in the lot-sizing literature \citep[i.e.][]{Kadipasa1997,Kilic2011}. A re-planning deployment proceeds as follows. The replenishment plans are made over the entire planning horizon, but only the imminent replenishment decision is implemented and a re-planning is done at the beginning of each period for the rest of the planning horizon. That is, at every period $n$, we determine a control policy for periods $\{n, n+1, \ldots, N\}$, but only implement the replenishment decision regarding period $n$. 

It is important to remark that policy parameters of the dynamic uncertainty strategy are independent of the inventory level on hand; i.e. the ($s$,$S$) levels of a non-stationary ($s$,$S$) policy remain optimal {\em regardless} of actual demand realisations. However, the same does not hold for static uncertainty and static-dynamic uncertainty strategies, for which the inventory position at the beginning of the planning horizon may affect the timing of replenishments that are scheduled in future periods, i.e. the $R$ parameter in the ($R$,$S$) and ($R$,$Q$) policies. This observation motivates the numerical study carried out in this section. Our aim is essentially to explore how the cost performance of different strategies is affected by the availability of realized demand information.

The re-planning approach applied in this study features similarities and differences from the conventional rolling horizon planning  \citep[p. 199]{Silver1998}. In rolling horizon planning a fixed length of time window is rolled forward after the imminent period's decision is implemented and the demand realization occurs. Rolling horizon planning for probabilistic demand is first implemented by \citet{Wemmerlov1984} and then revisited by \citet{Bookbinder1986} and \citet{Bookbinder1988}. These studies assume a long finite planning horizon so that fixed length of time window is rolled many times. As it is the case in the re-planning approach, towards the end of the horizon, these studies limit the length of time window with the number of remaining periods to capture end-of-horizon effects that are important for products featuring short life cycles. However, as the demand unfolds, these studies revise the demand forecasts for the rest of the planning horizon, in addition to updating the current inventory position. 

The practice of updating demand forecasts at the beginning of a given period --- i.e. modifying the distribution of demand in future periods --- is widespread among practitioners due to its cost effectiveness and practical relevance. However, this practice is justified if we assume that demand distributions in different periods are correlated. All the approaches contrasted in our study operate under the assumption that demands in different periods are independently distributed; under this assumption, demand forecast update is clearly not justified. Furthermore, it is worth remarking that by updating the distribution of demand in future periods one effectively changes the problem instance at hand. This strategy, although valuable in practice, is therefore not appropriate if one's aim is to reveal the added value of information about the realization of stochastic demands for a given problem instance. For this reason, in our study we simply update the current inventory position for re-planning purposes, leaving the demand forecasts unchanged. 

Since, as we remarked, policy parameters of the dynamic uncertainty strategy are independent of the inventory level on hand,\footnote{Note that this is not true if one updates demand forecasts at a given period, in which case new ($s$,$S$) levels will have to be recomputed for future periods since the problem instance has changed.}  in our study we only experiment on the re-planning implementations of static uncertainty and static-dynamic uncertainty strategies. 

In the following, we first present the experimental design adopted in our work, and then we discuss the results of the numerical study.

\subsection{Experimental setup}

We consider a planning horizon comprising 24 periods. Some of the studies under consideration do not immediately extend to generic demand distributions, e.g. \cite{Askin1981,Sox1997,Vargas2009} requires the cumulative demand distribution of all possible replenishment cycles, a convolution that often cannot be obtained in closed form; similarly, \cite{Bollapragada1999} only discuss the case of a Poisson demand and of a normal demand with a fixed coefficient of variation across periods; as the authors remark, their approach can be conceptually generalised to other distributions, however this requires the definition of tailor made search strategies \citep[][p. 578]{Bollapragada1999} and it is therefore not trivial. 

To enable a comparison across all approaches surveyed demand $d_n$ in each period $n = 1,\ldots,24$ is therefore assumed to be an independently normally distributed random variable with known expected value $\tilde{d}_n$ and standard deviation $\sigma_n = \rho\cdot\tilde{d}_n$, where $\rho$ denotes the coefficient of variation of the demand, which remains fixed over time as prescribed in \cite{Bollapragada1999}. As remarked, demand distributions are never updated, as this sort of action would be in contrast with the key assumption of demand independence over time periods, which is common across all methods we are working on. 

\begin{figure}[htbp]
\centering
\begin{minipage}[c]{0.49\linewidth}\includegraphics[width=\linewidth]{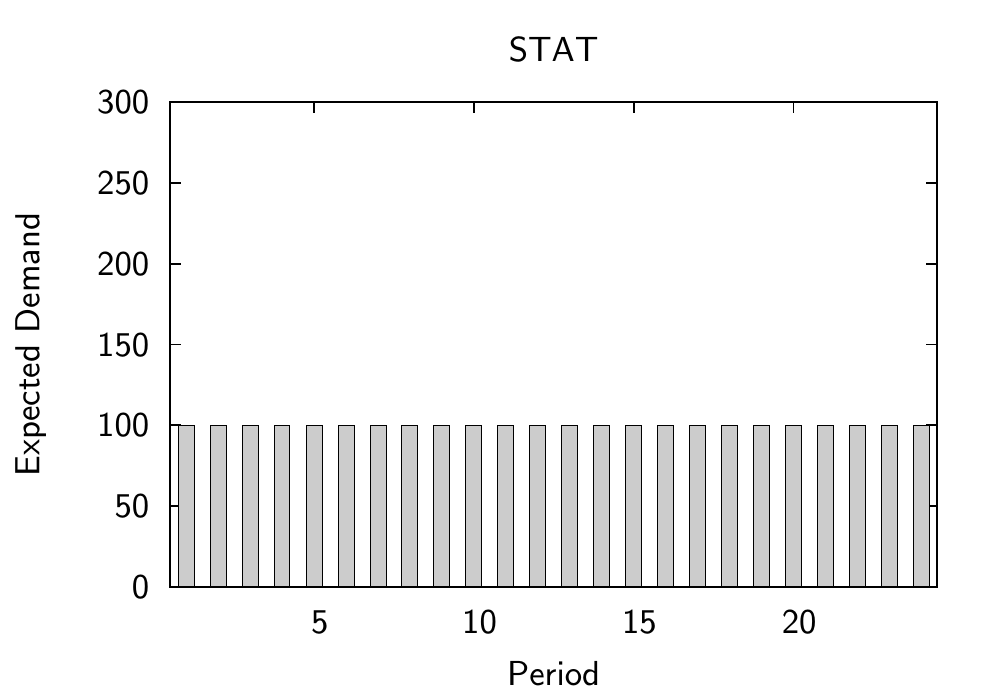}\end{minipage}
\begin{minipage}[c]{0.49\linewidth}\includegraphics[width=\linewidth]{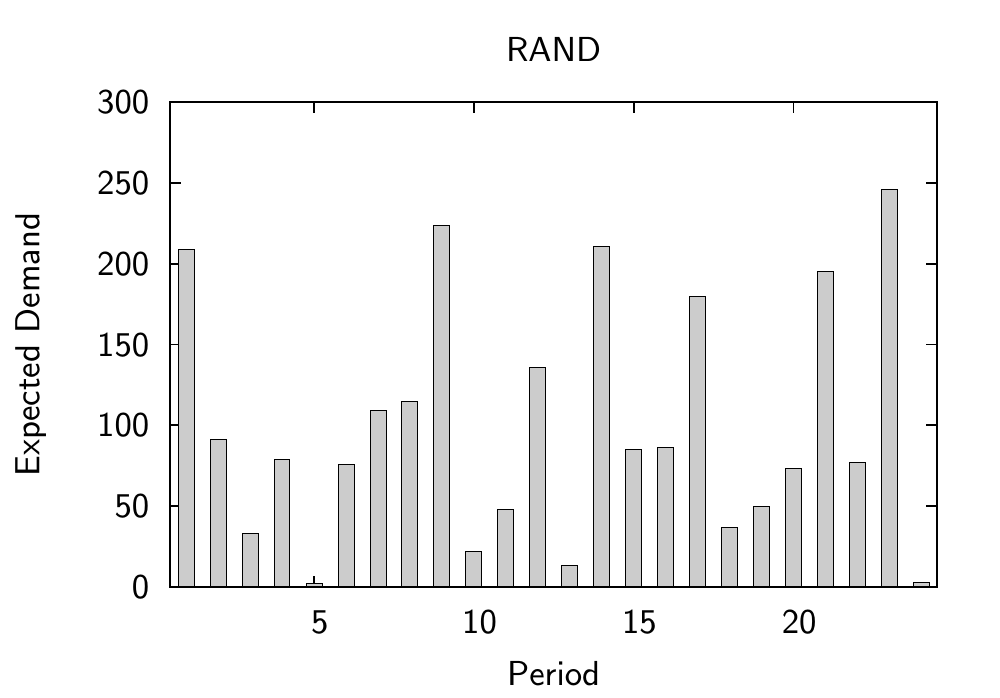}\end{minipage}\\
\begin{minipage}[c]{0.49\linewidth}\includegraphics[width=\linewidth]{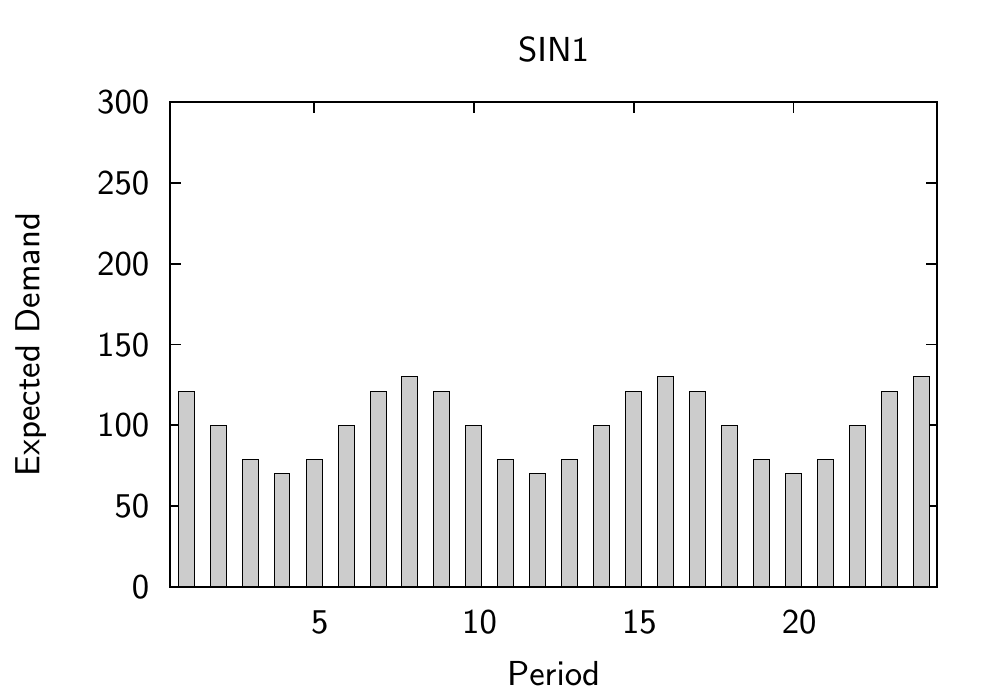}\end{minipage}
\begin{minipage}[c]{0.49\linewidth}\includegraphics[width=\linewidth]{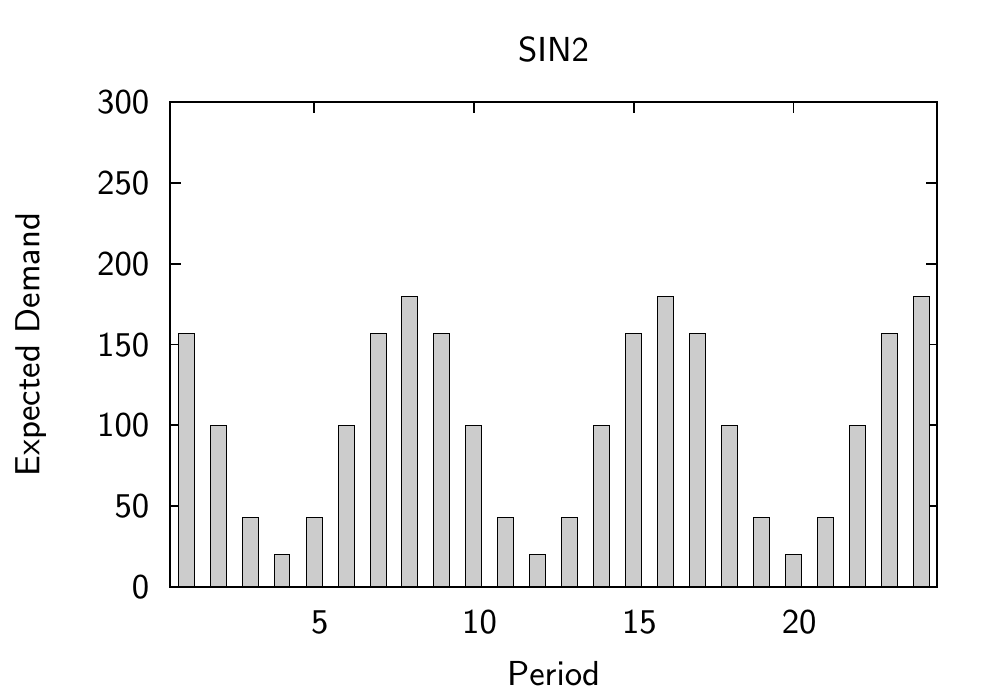}\end{minipage}\\
\begin{minipage}[c]{0.49\linewidth}\includegraphics[width=\linewidth]{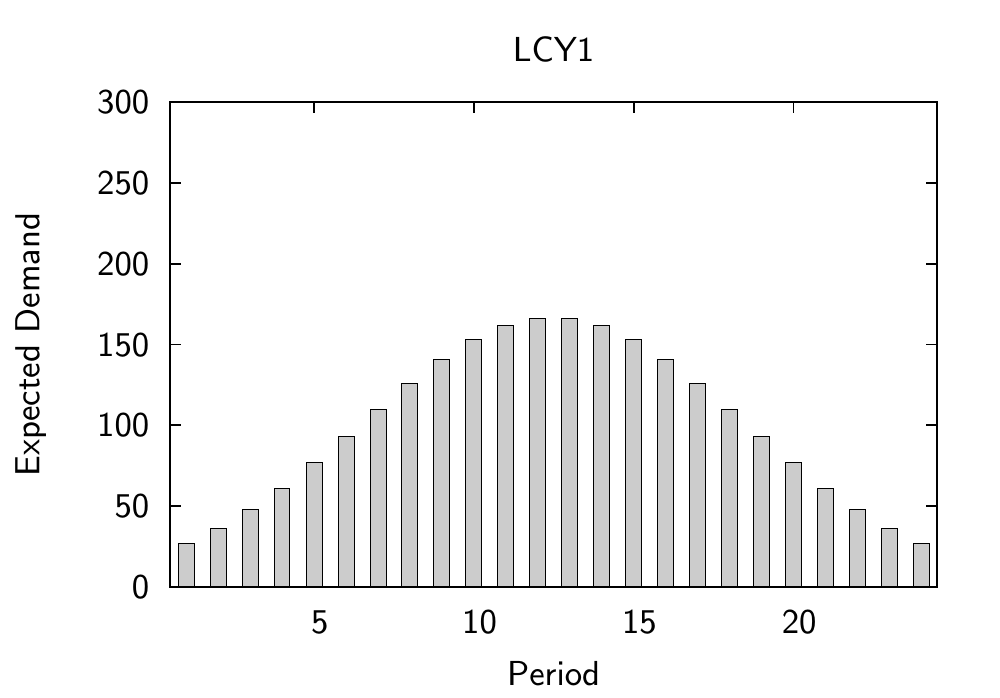}\end{minipage}
\begin{minipage}[c]{0.49\linewidth}\includegraphics[width=\linewidth]{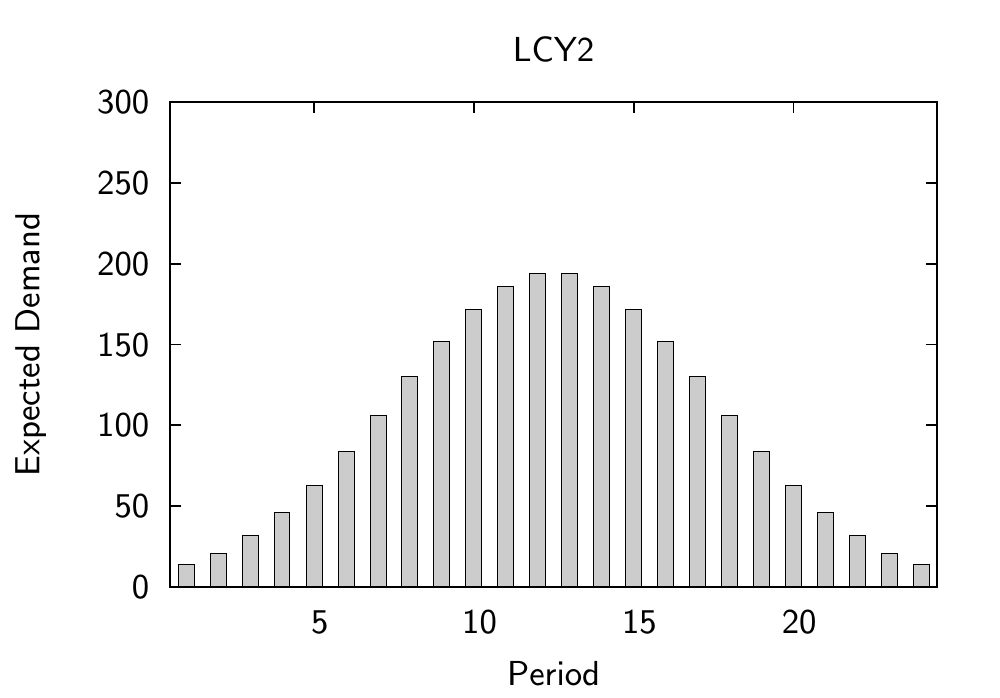}\end{minipage}
\caption{Demand patterns} \label{fig:demand_patterns}
\end{figure}

We consider 6 different patterns for the expected value of the demand in each period of the planning horizon: a stationary demand pattern (STA); an erratic pattern (RAND); two life cycle patterns, one with lower (LCY1) and one with higher (LCY2) variation of the expected demand over the planning horizon; and finally two sinusoidal patterns, one with weaker (SIN1) and one with stronger (SIN2) oscillations. We adapted these demand patterns from \citet{Berry1972}; the same patterns have been extensively used in the literature \citep{Tarim2006,citeulike:14063592,Rossi2008,Tunc2011}. Figure~\ref{fig:demand_patterns} graphically illustrates all demand patters analyzed in the numerical study. Note that the average demand per period is equal to 100 for all demand patterns. By varying the coefficient of variation $\rho\in\{0.10,0.20,0.30\}$, the fixed ordering cost $K\in\{250,500,1000,2000\}$, and the penalty cost $b\in\{2,5,10\}$ we generate a total of 216 test instances. 

We first obtain the optimal parameters of dynamic uncertainty strategy for each test instance by means of the stochastic dynamic programming approach of \citet{Scarf1960}. Then, we solve all instances by using each of the methods presented in Table~\ref{tab:overview}. For methods adopting static uncertainty and static-dynamic uncertainty strategy  \citep[i.e.][]{Vargas2009,Tarim2006,Rossi2012}, we also investigate the effect of adopting a re-planning approach with respect to the realized demand information. Finally, the heuristic proposed by \citet{Bollapragada1999} implicitly assumes an infinite planning horizon, and thus, it cannot account for the so-called end-of-horizon effect. In order to overcome this issue, the optimal policy parameters of dynamic uncertainty strategy are used for the last 8 periods of the planning horizon as it is done in their original study. 

We simulate the control policies obtained by implementing each of the aforementioned methods by using the common random numbers simulation strategy \citep[see e.g.][]{Sloan1990}. We implement a stopping rule so as to achieve an estimation error of $\pm$0.1\% of the expected total cost with 0.95 confidence probability. Average costs are then compared against the average cost of the optimal policy \citep{Scarf1960} obtained by the stochastic dynamic program  and the differences between the two are recorded as the percentage optimality gaps. 

We carry out sensitivity analysis by fixing one parameter at a time --- this is called the ``pivot'' parameter --- while varying others over their respective domain. The results of the numerical study are summarized in Table~\ref{table:summary} where the average optimality gap is reported for all problem instances characterized by the same pivot parameter. Also, in Figure~\ref{fig:boxplot_allresults} and in Figure~\ref{fig:boxplot_allresults_roll}, we provide a comprehensive set of boxplots for conventional and re-planning implementations of different methods. Note that the results regarding each method is given under the abbreviated name of the first author of the corresponding paper, i.e. Sox/Var, Tar, Ros, Ask, and Bol; and the suffix ``--R'' is added to those where a re-planning approach is adopted, i.e. Sox/Var--R, Tar--R, and Ros--R. Note that, \citet{Vargas2009} presents a special case implementation of \citet{Sox1997} by assuming a normally distributed demand. So, in the numerical setup of the current study these two are treated as a single method and abbreviated as Sox/Var.

\begin{sidewaystable}[htbp]
\centering
\scriptsize
\begin{tabular*}{1.04\textwidth}{@{\extracolsep{\fill}} lccccccccccccccccccccccc}
\toprule
\addlinespace
	& \multicolumn{2}{c}{Dynamic} & \multicolumn{8}{c}{Static-Dynamic}	& \multicolumn{2}{c}{Static}	\\
\cmidrule(r){2-3} \cmidrule(r){4-7} \cmidrule(r){8-9}
\addlinespace
	& Ask	& Bol	& Tar	& Tar--R	& Ros	& Ros--R	& Sox/Var	& Sox/Var--R \\
    \midrule
\multicolumn{3}{l}{\hspace{-10pt} Demand pattern ($\pi$)} \\
    STA   & 1.8   & 0.4   & 1.3   & 0.2   & 1.3   & 0.2   & 10.5  & 0.4 \\
    RAND  & 3.5   & 12.2  & 2.0   & 0.3   & 2.0   & 0.3   & 15.0  & 0.6 \\
    SIN1  & 2.8   & 1.9   & 1.5   & 0.2   & 1.4   & 0.2   & 10.7  & 0.4 \\
    SIN2  & 4.3   & 10.0  & 2.2   & 0.5   & 2.2   & 0.4   & 13.5  & 0.6 \\
    LCY1  & 4.6   & 1.9   & 1.2   & 0.1   & 1.1   & 0.1   & 12.8  & 0.4 \\
    LCY2  & 6.5   & 3.0   & 1.3   & 0.2   & 1.3   & 0.2   & 14.9  & 0.5 \\
\midrule
\multicolumn{3}{l}{\hspace{-10pt} Coefficient of variation ($\rho$)} \\
    0.1   & 3.6   & 6.0   & 0.2   & 0.1   & 0.2   & 0.0   & 5.2   & 0.1 \\
    0.2   & 4.0   & 4.8   & 1.2   & 0.2   & 1.2   & 0.2   & 12.4  & 0.4 \\
    0.3   & 4.2   & 3.9   & 3.3   & 0.5   & 3.3   & 0.4   & 21.2  & 0.9 \\
\midrule
\multicolumn{3}{l}{\hspace{-10pt} Setup cost ($K$)} \\
    250   & 3.3   & 4.3   & 2.0   & 0.4   & 2.0   & 0.4   & 25.0  & 0.7 \\
    500   & 3.1   & 4.3   & 1.8   & 0.3   & 1.8   & 0.3   & 14.3  & 0.5 \\
    1000  & 4.7   & 4.4   & 1.5   & 0.2   & 1.5   & 0.2   & 8.0   & 0.4 \\
    2000  & 4.5   & 6.6   & 0.9   & 0.1   & 0.9   & 0.1   & 4.3   & 0.3 \\
\midrule
\multicolumn{3}{l}{\hspace{-10pt} Penalty cost ($p$)} \\
    2     & 4.4   & 4.3   & 0.8   & 0.2   & 0.8   & 0.1   & 8.4   & 0.3 \\
    5     & 3.9   & 5.4   & 1.5   & 0.3   & 1.5   & 0.2   & 13.1  & 0.4 \\
    10    & 3.5   & 5.0   & 2.3   & 0.3   & 2.3   & 0.3   & 17.2  & 0.7 \\
\midrule
\multicolumn{3}{l}{\hspace{-10pt} All instances} \\
    AVG   & 3.9   & 4.9   & 1.6   & 0.3   & 1.6   & 0.2   & 12.9  & 0.5 \\
\bottomrule
\end{tabular*}
\caption{Average $\%$ optimality gaps of methods for different pivoting parameters}
\label{table:summary}
\end{sidewaystable}

\subsection{Discussion}

\begin{figure}[t!]
\centering
\subfigure[All implementations]{
\label{fig:boxplot_allresults}
\includegraphics[width=0.7\textwidth]{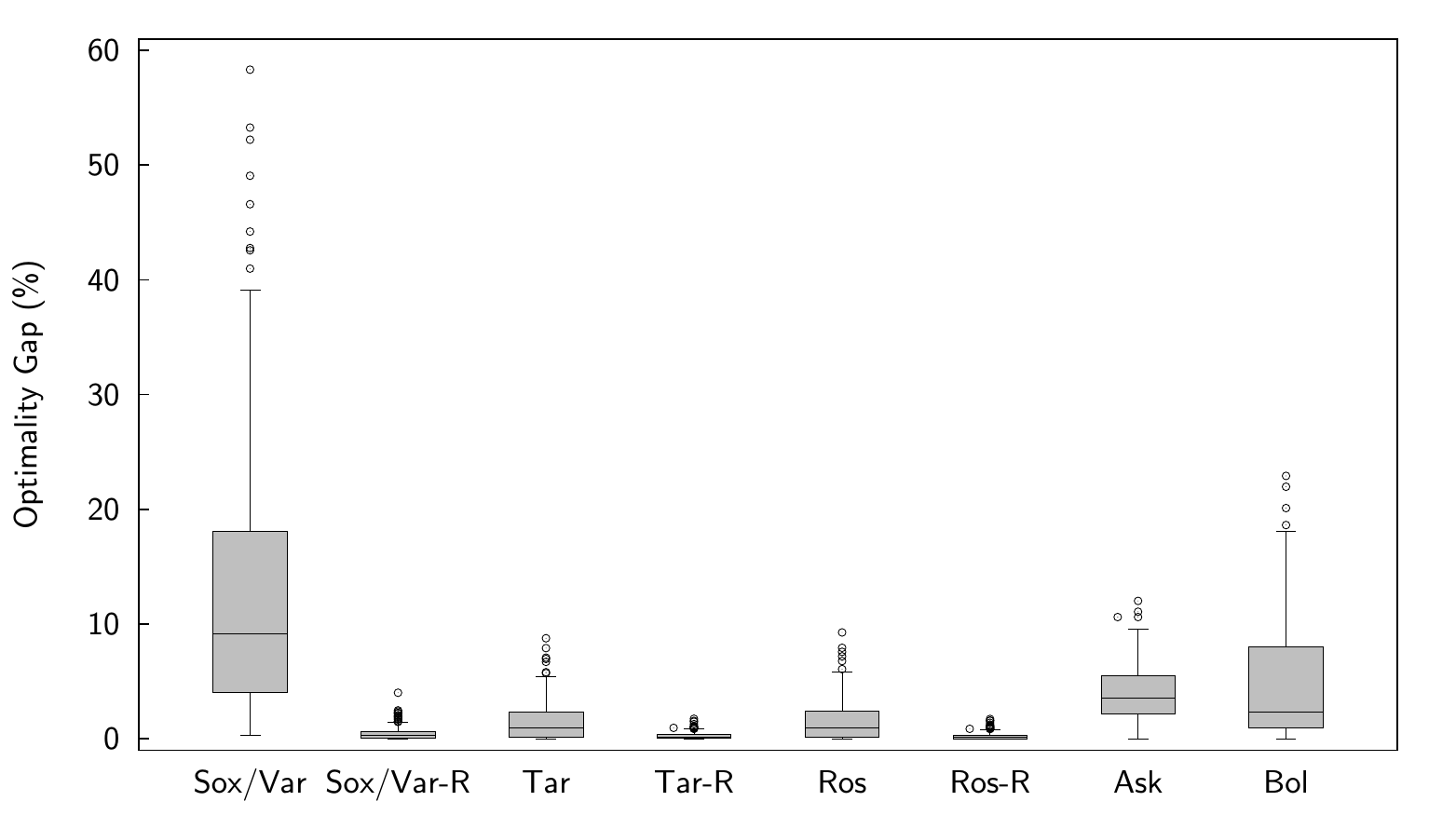}} 
~
\subfigure[Re-planning implementations -- A closer look]{
\label{fig:boxplot_allresults_roll}
\includegraphics[width=0.7\textwidth]{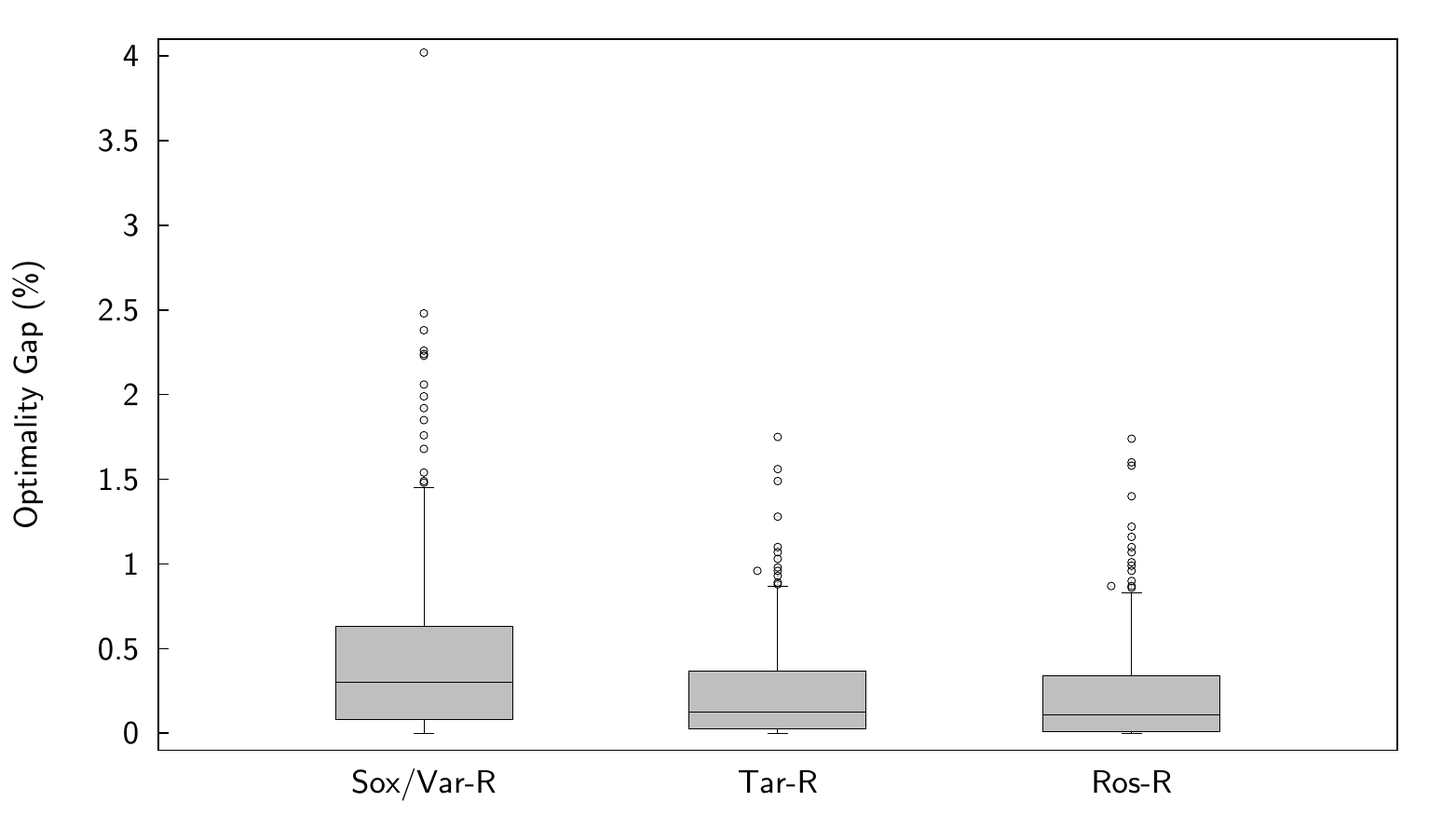}}
\caption{Boxplot comparison of methods considered}
\end{figure}

Below we first give a brief overview of the findings regarding each of the methods considered with and without re-planning. We draw some conclusions on the affects of varying demand and cost parameters. We refer to the results reported in Table~\ref{table:summary}, Figure~\ref{fig:boxplot_allresults} and Figure~\ref{fig:boxplot_allresults_roll}. 

The results indicate that the method of \citet{Sox1997} and \citet{Vargas2009} -- in the absence of a re-planning approach -- provides the worst cost performance. This result is rather consistent and it holds for all parameter settings yielding an average optimality gap of $13\%$. In Figure~\ref{fig:boxplot_allresults} we can see that the optimality gap may reach $60\%$ in the worst case. These results can mainly be attributed to the structure of the policy. Because both the timing and the size of replenishments are fixed at the beginning of the planning horizon, the static uncertainty strategy has no means to take recourse actions against demand uncertainty. Hence, uncertainty accumulates throughout the planning horizon and leads to a large optimality gap. 

It is possible to observe that the heuristics of \citet{Tarim2006} and \citet{Rossi2012} adopting the static-dynamic uncertainty strategy consistently yield the best cost performances among all methods considered, with an optimality gap around $2\%$. Furthermore, the deviation in the cost performance of these methods is rather low. This immediately shows that the static-dynamic uncertainty is very competitive despite the fact that it schedules replenishments in advance. 

\citet{Askin1981} and \citet{Bollapragada1999} yield optimality gaps around $4\%$ and $5\%$ respectively. These gaps are better than those produced by static uncertainty heuristics; however, they are much worse than those produced by static-dynamic uncertainty heuristics. This finding is a particularly interesting one. The dynamic uncertainty strategy is by definition a superior strategy as compared to the static-dynamic uncertainty strategy because it uses more information at the time of making replenishment decisions. However, this is true provided optimal policy parameters are available for both policies. From our study it appears that existing heuristics for the dynamic uncertainty policy do not produce parameters of sufficient quality to outpace the performance of state-of-the-art static-dynamic uncertainty heuristics.

We now concentrate on the effects of varying demand and cost parameters. The results suggest, in general, that all methods perform relatively better when demand pattern is rather steady. This is especially apparent when we look at variants of the same underlying pattern characterized by small and large oscillations. The average optimality gap of all methods significantly increase as we move from SIN1 and SIN2, or LCY1 to LCY2. This immediately shows that managing inventories is more difficult when demand is heavily non-stationary. 
 
We observe that performance of all methods deteriorates as demand variation increases. For instance, the optimality gap of the method of \citet{Sox1997} and \citet{Vargas2009} gradually increases from $5.15\%$ to $21.17\%$ as the coefficient of variation increases from $0.1$ to $0.3$. The only exception of this is the heuristic of \citet{Bollapragada1999}. This method performs relatively better in higher levels of demand uncertainty. This result can mainly be attributed to the fact that all methods except \citeauthor{Bollapragada1999}'s heuristic are adaptations of solution methods which were originally designed for the deterministic version of the lot sizing problem.  

It is not possible to make a general remark for all methods with respect to fixed ordering costs. Nevertheless, we see that the performance of the method of \citet{Sox1997} and \citet{Vargas2009}, as well as the heuristics of \citet{Tarim2006} and \citet{Rossi2012} improve when fixed ordering cost increases. This is not the case for dynamic uncertainty heuristics, i.e. \citet{Askin1981} and \citet{Bollapragada1999}. The effect of fixed ordering costs on these heuristics does not show a stable trend. 

We observe a similar picture when we look at the effect of stock-out penalty cost. The method of \citet{Sox1997} and \citet{Vargas2009} and the  heuristics of \citet{Tarim2006} and \citet{Rossi2012} perform worse for higher values of penalty cost. \citeauthor{Askin1981}'s (\citeyear{Askin1981}) heuristic, on the other hand, performs relatively better as the penalty cost increases. We cannot observe a stable trend on the performance of the heuristic of \citet{Bollapragada1999} with respect to penalty cost.

Having summarized the results obtained from conventional implementations of different methods, we turn our attention on their re-planning implementations. 
As remarked, we only target static-uncertainty and static-dynamic uncertainty strategies, since a re-planning approach is only of value if the policy parameters are dependent on the inventory level on hand. Hence, we experiment on the methods of \citet{Sox1997} and \citet{Vargas2009}, \citet{Tarim2006}, and \citet{Rossi2012}. We observe that the re-planning approach not only significantly improves the average cost performance of all these methods, but it also reduces the dispersion in the worst-case performance (see Figure~\ref{fig:boxplot_allresults_roll}). In particular, we see that the heuristics of \citet{Tarim2006} and \citet{Rossi2012} perform extremely well when deployed with a re-planning approach and both yield optimality gaps around $0.25\%$. On the other hand, probably the most remarkable result in our numerical experiments is that the re-planning implementation of the method of \citet{Sox1997} and \citet{Vargas2009} adopting static uncertainty strategy displays an excellent performance with an optimality gap around $0.5\%$ -- slightly higher than those of the heuristics based on static-dynamic uncertainty strategy. Note that this method displays an average optimality gap of $13\%$ without the re-planning deployment. This clearly shows that the re-planning approach makes it possible for the static uncertainty strategy to take suitable recourse actions against demand uncertainty. 

\section{Conclusions}
\label{sec:conclusions}

In this paper, we presented a cost-based performance evaluation of three lot-sizing strategies, i.e. static uncertainty, dynamic uncertainty, static-dynamic uncertainty, for non-stationary stochastic demand under the assumption of penalty cost and complete demand backordering. We compared six state-of-the-art heuristics \citep[i.e.][]{Scarf1960,Askin1981,Bollapragada1999,Tarim2006,Rossi2012, Sox1997, Vargas2009} on a common test bed and computed their optimality gap with respect to the optimal control policy obtained via stochastic dynamic programming. A similar cost comparison is carried out for a re-planning deployment of the static and static-dynamic heuristics. 

According to our numerical study the static-dynamic uncertainty strategy is the one that displays the closest cost performance to the optimal dynamic uncertainty strategy. In particular, heuristic methods of \citet{Tarim2006} and \citet{Rossi2012} yield an optimality gap around $2\%$ as opposed to the large optimality gap of static uncertainty strategy around $13\%$. Surprisingly, heuristic methods proposed for dynamic uncertainty strategy yield a cost performance in-between the two other strategies. Furthermore, their cost performance is very sensitive to the parameter settings and displays large dispersion against the parameter changes. These results suggest that the heuristics based on static-dynamic uncertainty strategy are strong alternatives for the optimal lot-sizing policy.

The numerical results tell that the static-dynamic uncertainty strategy is the one that displays the closest cost performance to the optimal dynamic uncertainty strategy. In particular, heuristic methods of \citet{Tarim2006} and \citet{Rossi2012} yield an optimality gap around $2\%$ as opposed to the large optimality gap of static uncertainty strategy around $13\%$. On the other hand, heuristic methods proposed for dynamic uncertainty strategy yield a cost performance in-between the two other strategies. However, their cost performance is very sensitive to the parameter settings and displays large dispersion against the parameter changes. 

The scene changes when re-planning with respect to realized demand information is incorporated. As one would expect, the re-planning approach improved the cost performance of all methods considered. Though, an outstanding impact is observed with the methods for static uncertainty strategy with an average optimality gap of $0.5\%$. This shows that when deployed with a re-planning approach the worst performing strategy, i.e. static uncertainty, turns into one of the strongest alternatives to the optimal dynamic uncertainty strategy.

\bibliographystyle{plainnat}
\bibliography{experiment}

\end{document}